\documentclass[10pt]{amsart}
\usepackage{amsfonts}
\usepackage{amsmath}
\usepackage{amsthm}
\usepackage{amssymb}
\usepackage{latexsym}
\usepackage{multicol}
\usepackage{verbatim}

\advance\textwidth by 1.2in \advance\oddsidemargin by -.6in
\advance\evensidemargin by -.6in
\newtheorem*{cor}{Corollary}
\newtheorem*{lem}{Lemma}
\newtheorem*{prop}{Proposition}

\theoremstyle{definition}

\theoremstyle{definition}
\newtheorem*{thm}{Theorem}

\newtheorem*{rem}{Remark}

\newenvironment{pf}{\proof}{\endproof}
\newcounter{cnt}
\newenvironment{enumerit}{\begin{list}{{\hfill\rm(\roman{cnt})\hfill}}{%
\settowidth{\labelwidth}{{\rm(iv)}}\leftmargin=\labelwidth%
\advance\leftmargin by
\labelsep\rightmargin=0pt\usecounter{cnt}}}{\end{list}}

\theoremstyle{remark}

\numberwithin{equation}{section} 

\def\wtl{{\rm wt}_\ell}
\def\wt{{\rm wt}}
\def\supp{{\rm supp}}
\def\chl{{\rm ch}_\ell}

\def\opl_#1^#2{\text{\tiny$\bigoplus\limits_{\text{\scriptsize$#1$}}^{\text{\scriptsize$#2$}}$}}
\def\hlie#1{\widehat{\mathfrak{#1}}}

\begin{document}

\newcommand{\thmref}[1]{Theorem~\ref{#1}}
\newcommand{\secref}[1]{Section~\ref{#1}}
\newcommand{\lemref}[1]{Lemma~\ref{#1}}
\newcommand{\propref}[1]{Proposition~\ref{#1}}
\newcommand{\corref}[1]{Corollary~\ref{#1}}
\newcommand{\remref}[1]{Remark~\ref{#1}}
\newcommand{\defref}[1]{Definition~\ref{#1}}
\newcommand{\er}[1]{(\ref{#1})}
\newcommand{\id}{\operatorname{id}}
\newcommand{\tensor}{\otimes}
\newcommand{\nc}{\newcommand}
\newcommand{\rnc}{\renewcommand}
\newcommand{\qbinom}[2]{\genfrac[]{0pt}0{#1}{#2}}
\nc{\cal}{\mathcal} \nc{\goth}{\mathfrak} \rnc{\bold}{\mathbf}
\renewcommand{\frak}{\mathfrak}
\newcommand{\desc}{\operatorname{desc}}
\newcommand{\Maj}{\operatorname{Maj}}
\renewcommand{\Bbb}{\mathbb}
\nc\bomega{{\mbox{\boldmath $\omega$}}}
 \nc\bvpi{{\mbox{\boldmath
$\varpi$}}}
 \nc\bpi{{\mbox{\boldmath
$\pi$}}}
 \nc\balpha{{\mbox{\boldmath
$\alpha$}}}

\newcommand{\lie}[1]{\mathfrak{#1}}
\makeatletter
\def\section{\def\@secnumfont{\mdseries}\@startsection{section}{1}%
  \z@{.7\linespacing\@plus\linespacing}{.5\linespacing}%
  {\normalfont\scshape\centering}}
\def\subsection{\def\@secnumfont{\bfseries}\@startsection{subsection}{2}%
  {\parindent}{.5\linespacing\@plus.7\linespacing}{-.5em}%
  {\normalfont\bfseries}}
\makeatother
\def\subl#1{\subsection{}\label{#1}}

\nc{\Cal}{\cal} \nc{\Xp}[1]{X^+(#1)} \nc{\Xm}[1]{X^-(#1)}
\nc{\on}{\operatorname} \nc{\ch}{\mbox{ch}} \nc{\Z}{{\bold Z}}
\nc{\J}{{\cal J}} \nc{\C}{{\bold C}} \nc{\Q}{{\bold Q}}
\renewcommand{\P}{{\cal P}}
\nc{\N}{{\Bbb N}} \nc\boa{\bold a} \nc\bob{\bold b} \nc\boc{\bold
c} \nc\bod{\bold d} \nc\boe{\bold e} \nc\bof{\bold f}
\nc\bog{\bold g} \nc\boh{\bold h} \nc\boi{\bold i} \nc\boj{\bold
j} \nc\bok{\bold k} \nc\bol{\bold l} \nc\bom{\bold m}
\nc\bon{\bold n} \nc\boo{\bold o} \nc\bop{\bold p} \nc\boq{\bold
q} \nc\bor{\bold r} \nc\bos{\bold s} \nc\bou{\bold u}
\nc\bov{\bold v} \nc\bow{\bold w} \nc\boz{\bold z}

\nc\ba{\bold A} \nc\bb{\bold B} \nc\bc{\bold C} \nc\bd{\bold D}
\nc\be{\bold E} \nc\bg{\bold G} \nc\bh{\bold H} \nc\bi{\bold I}
\nc\bj{\bold J} \nc\bk{\bold K} \nc\bl{\bold L} \nc\bm{\bold M}
\nc\bn{\bold N} \nc\bo{\bold O} \nc\bp{\bold P} \nc\bq{\bold Q}
\nc\br{\bold R} \nc\bs{\bold S} \nc\bt{\bold T} \nc\bu{\bold U}
\nc\bv{\bold V} \nc\bw{\bold W} \nc\bz{\bold Z} \nc\bx{\bold X}
\title[]{Characters of fundamental representations of quantum affine algebras.}

\author{Vyjayanthi Chari and Adriano A. Moura }
\address{Department of Mathematics, University of
California, Riverside, CA 92521.} \email{chari@math.ucr.edu,
adrianoam@math.ucr.edu}
\begin{abstract} We give  closed formulae for the
$q$--characters of the fundamental representations of the quantum
loop algebra of a classical Lie algebra, in terms of a family of
partitions satisfying some simple properties. We also give the
multiplicities of the eigenvalues of the imaginary subalgebra in
terms of these partitions.

\end{abstract}

\maketitle \setcounter{section}{0}

\section*{Introduction}

In this paper we study the $q$--characters of the fundamental
finite--dimensional  representations of the  quantum loop algebra
$\bu_q$ associated to a classical simple Lie algebra. The notion
of $q$--characters defined in \cite{FR}  is analogous to the usual
notion of a character of a finite--dimensional representation of a
simple Lie algebra. These characters and their generalizations
have been studied extensively  \cite{FM}, \cite{He}, \cite{Na}
using combinatorial and geometric methods. A more representation
theoretic approach  was developed in \cite{CM}. In particular,
that paper approached the problem of studying whether the
$q$--characters admitted a Weyl group invariance which was
analogous to the invariance of    characters of
finite--dimensional representations of  simple Lie algebras. In
the quantum case, it is reasonable to expect that the Weyl group
be replaced by the braid group, \cite{BP}, \cite{C1}, \cite{FR}
but it is easy to see that this is false even for $sl_2$.
 However, it was shown in \cite{CM} that in a suitably modified way,
 the $q$--characters of the
fundamental representations of the quantum loop algebra of a
classical Lie algebra do admit an  invariance under the braid
group action. It was also shown that the $q$--character of such
representations could then be calculated in a certain inductive
way.

In this paper, we use that inductive method to give closed
formulas for the $q$--characters of all the fundamental
representations of the quantum loop algebras of a classical simple
Lie algebra.  To describe the results a bit further, recall that
the quantum loop algebra admits a commutative subalgebra
$\bu_q(0)$ corresponding to the imaginary root vectors. Any
finite--dimensional representation $V$ of the quantum  loop
algebra, breaks up as a direct sum of
 generalized eigenspaces  for the action of $\bu_q(0)$.
These are called the $\ell$--weight spaces and the eigenvalues
corresponding to the non--zero eigenspaces are called  the
$\ell$--weights of the representation. The $\ell$--weights  lie in
a free abelian multiplicative group  $\cal{P}_q$. Let  $\bz[\cal
P_q]$ be the integral group ring over $\cal P_q$ and for
$\bvpi\in\cal P_q$, let $V_\bvpi$ be the corresponding eigenspace
of $V$. The element of $\bz[\cal P_q]$ defined by, $${\rm
ch}_\ell(V)=\sum_{\bvpi\in\cal{P}_q}\dim(V_\bvpi)\ e(\bvpi),$$
 is called the $q$--character of $V$. If  $P$ is  the usual weight
 lattice of the simple Lie algebra, then it was shown in \cite{CM}
 that there exists a canonical group homomorphism $\wt:\cal P_q\to
 \cal P$.

 Assume now that $V$ is a fundamental representation of the
 quantum loop algebra. Roughly speaking, this means that $V$
 corresponds to a canonical  generator of $\cal P_q$. It was shown in
 \cite{CM} that the problem of determining the $\ell$--weights of
 $V$ is reduced to determining $V_\bvpi$ where $\wt(\bvpi)$ is in
 the dominant chamber $P^+$ of $P$. Assume from now on that
 $\wt(\bvpi)\in P^+$. We give explicit formulas for $\bvpi$ with
 $V_\bvpi\ne 0$. In the case of $B_n$, $C_n$, we see as a
 consequence that  $\dim V_\bvpi = 1$ (this was proved by
 different methods in \cite{He}). In the case of $D_n$ it can
 happen that $\dim(V_\bvpi)>1$ and we compute this dimension in
 Section 5. The idea is to show that every $\ell$--weight $\bvpi$ comes from  a partition
 $\boj$
with certain properties and we find that  $\dim
V_\bvpi=2^{M_\boj}$ where $M_\boj$ is defined in a canonical way
in Section 5.

\section{Preliminaries}

\subsection{} Let $\frak g$ be a complex finite--dimensional simple
Lie algebra  of rank
 $n$ and let $\frak h$ be a Cartan subalgebra of $\lie g$.
 Set $I=\{1,2,\cdots ,n\}$ and  let $\{\alpha_i:i\in I\}$ (resp.  $\{\omega_i:i\in I\}$)
  be the set of simple roots
 (resp.  fundamental weights) of $\frak g$ with respect to $\frak h$. Let also $\check{\alpha}_i$ denote the simple co-roots. As usual,
 $Q$,
(resp. $P$) denotes the  root (resp. weight) lattice of $\frak g$,
$Q^+ = \sum_{i=1}^n \bn \alpha_i$, and $P^+= \sum_{i=1}^n \bn \omega_i$.
 Let  $W$ be  the Weyl group of $\lie g$ generated by
simple reflections $\{s_i:i\in I\}$. For $w\in W$, let $\ell(w)$
denote the length of a reduced expression for $w$. Given
$\lambda=\sum_{i\in I} \lambda_i\omega_i\in P^+$ let $W(\lambda)$
be the subgroup of $W$ generated by $\{s_i: i\in I,\
\lambda_i=0\}$ and let   $W_\lambda$ be the set of left coset
representatives of $W/W(\lambda)$ of minimal length. The braid
group $\cal B$  associated to $\lie g$ is generated by elements
$T_i$, $i\in I$ and relations \begin{align*} T_iT_j &=T_jT_i,\ \
\text{if}\ \ a_{ij} =0,\\ T_iT_jT_i& =T_jT_iT_j,\ \ \text{if}\ \
a_{ij}a_{ji} =1,\\ (T_iT_j)^2&= (T_jT_i)^2,\ \ \text{if}\ \
a_{ij}a_{ji}=2,\\ (T_iT_j)^3&= (T_jT_i)^3,\ \ \text{if}\  \
a_{ij}a_{ji} =3,\end{align*} where $i,j\in\{1,2,\cdots ,n\}$ and
$A=(a_{ij})$ $1\le i,j\le n$ is the Cartan matrix of $\lie g$. For
$i\in I$, fix integers $d_i\in\bn $ minimal   such that
$d_ia_{ij}=d_ja_{ji}$ for all $j\in I$. Given $w\in W$ and a
reduced expression $w=s_{i_1}\cdots s_{i_k}$ let
$T_w=T_{i_1}\cdots T_{i_k}$ be the corresponding element of $\cal
B$. It is well--known that $T_w$ is independent of the choice of
the reduced expression.

\subsection{} Let $q\in\bc^\times$ and assume that $q$ is not a root
of unity. For $r,m\in\bn$, $m\ge r$, define complex numbers,
\begin{equation*}
[m]_q=\frac{q^m -q^{-m}}{q -q^{-1}},\ \ \ \ [m]_q!
=[m]_q[m-1]_q\ldots [2]_q[1]_q,\ \ \ \ \left[\begin{matrix} m\\
r\end{matrix}\right]_q = \frac{[m]_q!}{[r]_q![m-r]_q!}.
\end{equation*}
 Set $q_i=q^{d_i}$ and $[m]_i=[m]_{q_i}$.

 Let $\cal P_q$ be the (multiplicative) subgroup of $\bc(u)^n$ generated by the elements,
  $\bomega_{i,a}$, $i\in I$, $a\in\bc^\times$, where
$\bomega_{i,a}$ is the $n$--tuple of elements  in $\bc(u)$ whose
$i^{th}$ entry is $1-au$ and all other entries 1. The elements
$\bomega_{i,a}$ are called  $\ell$--fundamental weights. It is
obvious that $\cal{P}_q$ is generated freely as an abelian group
by the fundamental $\ell$--weights. $\cal P_q$ is called the
$\ell$--weight lattice. Given any element $\bvpi\in\cal P_q$ and
$1\le j\le n$, let $\bvpi_j$ be the $j^{th}$ entry of $\bvpi$.

 Let $\cal{P}_q^+$ be the monoid
generated by $1$ and the elements $\bomega_{i,a}$, $i\in I$,
$a\in\bc^\times$, clearly $\cal{P}_q^+$ consists  of $n$--tuples
of polynomials with constant term one and  an element of $\cal
P_q^+$ is called an  $\ell$--dominant weight.  Let
$\wt:\cal{P}_q\to P$ be the group homomorphism defined by
extending, $\wt(\bomega_{i,a})=\omega_i.$

The group $\cal B$ acts on
 on $ \cal{P}_q$ as follows \cite{BP},\cite{C1},\cite{FR}: for $i\in I$ and $\bvpi=(\varpi_1,\cdots ,\varpi_n)
\in\cal{P}_q$, we have
\begin{align*}
(T_i\bvpi)_j & =\varpi_j,\ \  {\text{if}}\ a_{ji}=0,\\
(T_i\bvpi)_j & =\varpi_j(u)\varpi_i(q_iu),\ \  {\text{if}}\
a_{ji}=-1,\\ (T_i\bvpi)_j & =\varpi_j(u)\varpi_i(q^3u)
\varpi_i(qu),\ \ {\text{if}}\ a_{ji}=-2,\\ (T_i\bvpi)_j &=
\varpi_j(u)\varpi_i(q^5u)\varpi_i(q^3u)\varpi_i(qu),\ \
{\text{if}}\ a_{ji}=-3,\\
 (T_i\bvpi)_i&
 =\frac{1}{\varpi_i(q_i^2u)}.\end{align*}

For $i\in I$, set $$\balpha_{i,a}= (T_{i}(\bomega_{i,
a}))^{-1}\bomega_{i, a}, $$  and let $\cal{Q}_q$ be the subgroup
of $\cal{P}_q$ generated by the $\balpha_{i,a}$. Let $\cal{Q}_q^+$
the monoid generated by $1$ and $\balpha_{i,a}$, $i\in I$,
$a\in\bc^\times$, and $\cal Q_q^-=(\cal Q_q^+)^{-1}$.

\subl{} The quantum loop algebra  $\bu_q$ of $\lie g$ is  the
algebra with generators $x_{i,r}^{{}\pm{}}$ ($i\in I$, $r\in\bz$),
$K_i^{{}\pm 1}$ ($i\in I$), $h_{i,r}$ ($i\in I$, $r\in
\bz\backslash\{0\}$) and the following defining relations:
\begin{align*}
   K_iK_i^{-1} = K_i^{-1}K_i& =1, \ \
 K_iK_j =K_jK_i,\\  K_ih_{j,r}& =h_{j,r}K_i,\\
 K_ix_{j,r}^\pm K_i^{-1} &= q_i^{{}\pm
    a_{ij}}x_{j,r}^{{}\pm{}},\ \ \\
  [h_{i,r},h_{j,s}]=0,\; \; & [h_{i,r} , x_{j,s}^{{}\pm{}}] =
  \pm\frac1r[ra_{ij}]_{q^i}x_{j,r+s}^{{}\pm{}},\\
 x_{i,r+1}^{{}\pm{}}x_{j,s}^{{}\pm{}} -q_i^{{}\pm
    a_{ij}}x_{j,s}^{{}\pm{}}x_{i,r+1}^{{}\pm{}} &=q_i^{{}\pm
    a_{ij}}x_{i,r}^{{}\pm{}}x_{j,s+1}^{{}\pm{}}
  -x_{j,s+1}^{{}\pm{}}x_{i,r}^{{}\pm{}},\\ [x_{i,r}^+ ,
  x_{j,s}^-]=\delta_{i,j} & \frac{ \psi_{i,r+s}^+ -
    \psi_{i,r+s}^-}{q_i - q_i^{-1}},\\
\sum_{\pi\in\Sigma_m}\sum_{k=0}^m(-1)^k\left[\begin{matrix}m\\k\end{matrix}
\right]_{i}
  x_{i, r_{\pi(1)}}^{{}\pm{}}\ldots x_{i,r_{\pi(k)}}^{{}\pm{}} &
  x_{j,s}^{{}\pm{}} x_{i, r_{\pi(k+1)}}^{{}\pm{}}\ldots
  x_{i,r_{\pi(m)}}^{{}\pm{}} =0,\ \ \text{if $i\ne j$},
\end{align*}
for all sequences of integers $r_1,\ldots, r_m$, where $m
=1-a_{ij}$, $\Sigma_m$ is the symmetric group on $m$ letters, and
the $\psi_{i,r}^{{}\pm{}}$ are determined by equating powers of
$u$ in the formal power series $$\sum_{r=0}^{\infty}\psi_{i,\pm
r}^{{}\pm{}}u^{{}\pm r} = K_i^{{}\pm 1}
{\text{exp}}\left(\pm(q_i-q_i^{-1})\sum_{s=1}^{\infty}h_{i,\pm s}
u^{{}\pm s}\right).$$ Let $\bu_q(\lie g)$ be the subalgebra of
$\bu_q$ generated by the elements $x_{i,0}^\pm$,  $K_i^{\pm 1}$
for $1\le i\le n$.

 \vskip 12pt

\subsection{}
For $i\in I$,  set
\begin{equation*}
h^\pm_i(u)=\sum_{k=1}^\infty \frac{q^{\pm k}h_{i,\pm
k}}{[k]_{i}}u^k,\end{equation*} and define elements $P_{i,\pm k}$,
$i\in I$, $k\in\bz$, $k\ge 0$, by the generating series,
\begin{equation}\label{hp}
P^\pm_i(u)=\sum_{k=0}^\infty P_{i,\pm
k}u^k=\exp(-h^\pm_i(u)).\end{equation}

Let $\bu_q^\pm(0)$ be the subalgebra of $\bu_q$ generated by the
elements $h_{i,\pm k}$ $i\in I$, $k\in\bz$, $k>0$, or equivalently,
the subalgebra generated by the elements $P_{i,\pm k}$, $i\in I$,
$k\in\bz$, $k>0$, and let $\bu_q(0)$ be the subalgebra generated by
$\bu_q^\pm(0)$.  An element $\bvpi=(\varpi_1,\cdots
,\varpi_n)\in\cal P_q$ can be regarded as an element of
$\text{Hom}(\bu_q(0),\bc)$ by extending the assignment, $$
\bvpi(P_i^\pm(u))= \varpi_i^\pm(u) ,$$ where
$\varpi_i^+(u)=\varpi_i$, $\varpi_i^-=u^{\text{deg}
\varpi_i}\varpi_i(u^{-1})/ (u^{\text{deg}
\varpi_i}\varpi_i(u^{-1}))|_{u=0}$.

\subl{}

 Given a $\bu_q$-module $V$ and $\mu=\sum_i\mu_i\omega_i\in
P$, set
\begin{equation*} V_\mu=\{ v\in V: K_i.v =q_i^{\mu_i}v ,\ \
\forall \ i\in I\}.\end{equation*} We say that $V$ is a module of
type 1 if
\begin{equation*} V=\bigoplus_{\mu\in P}V_\mu.\end{equation*}
Set $$\wt(V)= \{\mu\in P: V_\mu\ne 0\}, $$ and given $v\in V_\mu$
set $\wt(v)=\mu$.
An element
 $\bvpi\in\cal{P}_q$
is an $\ell$--weight of $V$ if there exists a non--zero element
$v\in V$ such that $$(P_{i, \pm r}-(\varpi_i^\pm)_r)^N v=0, \ \ \
N\equiv N(i,r,v)\in\bz^+,$$ for all $i\in I$ and $r\in\bz^+$ and
$v$ is called  an $\ell$--weight vector in $V$ with $\ell$--weight
$\bvpi$. Let $V_\bvpi$ be  the subspace of $V$ spanned by
$\ell$--weight vectors with $\ell$--weight $\bvpi$. If $V$ is a
finite--dimensional $\bu_q$--module, then, $$V=\opl_{\bvpi\in\cal
P_q}^{}V_{\bvpi},\ \ \ \ V_\mu=\opl_{\bvpi\in\cal P_q}^{}V_{\bvpi}
\cap V_\mu .$$
 Denote by $\wtl(V)$ the set of $\ell$--weights of $V$ and
define $\wtl(v)$ in the obvious way.

\subl{} Let $\cal C_q$ be the category of finite-dimensional
$\bu_q$--modules of type 1. A module $V\in\cal{C}_q$  is
$\ell$--highest weight with $\ell$--highest weight
$\bvpi\in\cal{P}_q$  if there exists a non--zero vector $0\ne v\in
V$ such that $V=\bu_qv$ and,
\begin{equation}\label{fdrel} x_{i,r}^+v =0,\ \
P^\pm_i(u)v=(\bvpi)^\pm_iv,\ \ K_i^{\pm
1}v=q^{\pm\wt\bvpi(\check{\alpha}_i)} v,\ \
(x_{i,r}^-)^{\wt\bvpi(\check{\alpha}_i)+1} v=0,\end{equation} for
all $i\in I$, $r\in\bz$. The element $v$ is called the
$\ell$--highest weight vector.

Any $\ell$--highest weight module has a unique irreducible
quotient which is also a highest weight module with the same
highest weight.  There exists a bijective correspondence  between
elements of $\cal{P}_q^+$ and isomorphism classes of  irreducible
finite--dimensional modules, \cite{CPbanff}. Given
$\bomega\in\cal{P}_q^+$, let $V(\bomega)\in\cal{C}_q$ be an
element in the corresponding isomorphism class, and let
$v_\bomega$ be the $\ell$--highest weight vector. Then,
$V(\bomega)_{\wt\bomega}=\bc v_\bomega.$

\subl{} From now on we suppose $\lie g$ is of classical  type, i.e
$\lie g$ is of type $A_n$, $B_n$, $C_n$ or $D_n$. The following
result was proved in \cite{CM}.
\begin{thm}\label{wconjfund}  Let $i\in I$, $a\in\bc^\times$ and $V=V(\bomega_{i,a})$. Assume that
$\bvpi\in\wt_\ell(V)$ is such that
$\wt(\bvpi)=\lambda\in P^+$.
\begin{enumerit}

\item[(i)] For all $w\in W_\lambda$ we have
$$\dim(V_\bvpi)=\dim(V_{T_w\bvpi}),$$
and
$$T_w(\wt_\ell(V_\lambda))=\wt_\ell(V_{w\lambda}).$$
\item[(ii)] Suppose that $\bvpi\ne \bomega_{i,a}$.  There exists
$\bvpi'\in\wt_\ell(V)$, $\mu=\wt(\bvpi')\in P^+$,  $w\in W_\mu$,
$j\in I$ with $\ell(s_jw)=\ell(w)+1$, and $c\in \bc^\times$ such
that $$(T_w(\bvpi'))_j = (1-cu)(1-c'u) \qquad\text{and}\qquad
\bvpi=T_w(\bvpi')(\balpha_{j,c})^{-1}$$ for some $c'\ne cq_j^{2}$,
and    $$\dim(V_{\bvpi})\ge 2 \qquad \text{if}\qquad c=c'.$$
Further, for all  $v\in V_{T_w(\bvpi')}$ and  $s\in\bz$,
\begin{equation}\label{multp} x_{j,s}^-v \in V_{\bvpi} +
V_{T_w(\bvpi')(\balpha_{j,c'})^{-1}}.\end{equation}
\end{enumerit}\hfill\qedsymbol
\end{thm}

\begin{cor} We have
$$\chl(V)=\sum_{\lambda\in P^+}\ \sum_{w\in W_{\lambda}} \
\sum_{\bvpi\in \wtl(V_\lambda)} (\dim V_{\bvpi})
e(T_w(\bvpi)).$$\hfill\qedsymbol
\end{cor}
>From now on, we will let  $\bvpi$ also denote the element
$e(\bvpi)$ of $\bz[\cal P_q]$. Notice that since the group $\cal
P_q$ is multiplicative, this should cause no confusion.

\subl{} It follows from the corollary that if $V(\bomega_{i,a})$
is a minuscule representation of $\lie g$, i.e
$V(\bomega_{i,a})_\lambda=0$ for all $\lambda\in P^+$ with
$\lambda< \omega_i$, then $$\chl(V(\bomega_{i,a}))= \sum_{w\in
W_{\omega_i}}\ T_w(\bomega_{i,a}).$$ It follows from
\cite{CPbanff} that this is the case for all fundamental
representations of $A_n$, the spin nodes for the orthogonal
algebras, and the natural representations of $C_n$ and $D_n$. In
the rest of the paper we consider the remaining cases.

\subl{} We conclude this section with a stronger version of
Theorem \ref{wconjfund}(ii). Let $\bu_{q_j}(\hlie g_j)$ be the
subalgebra of $\bu_q$ generated by the elements $x_{j,m}^\pm$,
$h_{j,s}$, $K_j^{\pm 1}$, $m,s\in\bz$, $s\ne 0$. It is known that
$\bu_{q_j}(\hlie g_j)$ is isomorphic to $\bu_{q_j}(\hlie{sl}_2)$.
\begin{prop}\label{dim}  Let $i\in I$, $a\in\bc^\times$,
$V=V(\bomega_{i,a})$. Let $\bvpi'\in\wtl(V)$ satisfy the
following: $\mu=\wt(\bvpi')\in P^+$,  $$(T_w(\bvpi'))_j =
(1-cu)(1-c'u), $$ for some $c'\ne cq_j^{2}$ and  $w\in W_{\mu}$,
$j\in I$ satisfying $(w(\mu)-\alpha_j)\in P^+$. Set,
 $\bvpi=T_w(\bvpi')(\balpha_{j,c})^{-1}.$ Then
$$\dim (V_{\bvpi})\ge \dim (V_{\bvpi'}).$$ Moreover, if $c=c'$,
then $$\dim (V_{\bvpi})\ge 2\dim (V_{\bvpi'}).$$
\end{prop}

\begin{pf} First note that if $\mu=\wt(\bvpi')\in P^+$ then
$\mu=\omega_r$ for some $r\le i$, (see \cite[Section 1]{CM} for
instance). Further, since $w\in W$ and $j\in I$ are such that
$\ell(s_jw)=\ell(w)+1$ it follows that
$w\omega_r+\alpha_j\notin\wt(V)$. Since
$\wt(T_w(\bvpi'))=w\omega_r$, it follows that
$$x_{j,m}^+V_{T_w(\bvpi')}=0,\ \ \forall\ m\in\bz.$$ Choose a
basis $\{v_1, \cdots, v_p\}$ of $V_{T_w(\bvpi')}$ such that
$$\bu_q(0)v_m\in\sum_{s\le m}\bc(q)v_s.$$ Let $U_m
=\bu_{q_j}(\hlie g_j)v_m$. Then $U_m/U_{m-1}$ is an
$\ell$--highest weight module for $\bu_{q_j}(\hlie g_j)v_m$ with
highest weight $(1-cu)(1-c'u)$ and hence  by \cite{CM},  \cite{FR}
there exists a unique (up to scalar multiple) non--zero element in
the span of $\{x_{j,r}^-v_m:r\in\bz\}$ which is an eigenvector
(modulo $U_{m-1}$) for the $P_{j,r}$ with eigenvalue $\varpi_j$ if
$c\ne c'$ and two linearly independent elements if $c=c'$.
Equation \eqref{multp} of  Theorem \ref{wconjfund}(ii) now implies
that such vectors are in $V_\bvpi$ and hence the proposition is
proved.
\end{pf}

\begin{rem}
It will actually follow from  Theorem \ref{main} and its proof
that equality holds in Proposition \ref{dim}.
\end{rem}

\section{Closed Formulae for $q$--characters}

In this section we state the main theorem which gives closed
formulas for the $\ell$--weights $\bvpi$ with $\wt(\bvpi)\in P^+$
of the fundamental representations of quantum affine algebras.

\subl{} Assume that the Dynkin diagram of $\lie g$ is labeled as
in \cite{Bo}. Throughout this section we shall assume that we have
fixed an integer $i$ such that
\begin{align*}
1< i\le n &&\text{ if } \lie g = C_n,\\
 1\le  i<n && \text{ if } \lie g = B_n,\\
 1<i\le n-2 && \text{ if } \lie g = D_n.\end{align*}

Define a subset $I_i$ of $I$ by,
\begin{align}\label{iib}
I_i =&\ \{r: 0\le r\le i\},&&  \text{ if } \lie g = B_n,\\
\label{iicd}=&\ \{r: 0\le r\le i, r\equiv i \ \rm {mod}\  2\},&&
\text{ if } \lie g = C_n,D_n.
\end{align}

>From now on, given $r\in I_i$, we shall denote by $M$ the greatest
integer less than or equal to $(i-r)/2$.

\subl{} For  $r\in I_i$, let $\bj_{k,r}$  be the set of partitions
$r<j_1<j_2<\cdots <j_k\le n$ of length $k$ and satisfying
\begin{align*} j_s&\le n-i+r+2s-1, \ \ 1\le s\le k, &&
\text{ if } \lie g = C_n,\\ j_k&<n, && \text{ if } \lie g =
D_n.\end{align*} Set
\begin{align}
\label{jc} \bj_r  &=\bj_{M,r} && \text{ if } \lie g = C_n,\\
\label{jbd} &= \cup_{0\le k\le M} \bj_{k,r} &&\text{ if } \lie g =
B_n, D_n,
\end{align}
where $\bj_{0,r}$ consists of the empty partition.

\subl{}  Given $j\in I$, $r\in I_i$ and an integer $2d_1s\in\bn$,
define elements $\bpi_r(j,s)\in\cal P_q$ by
\begin{equation*}
\bpi_r(j,s) =  \bomega_{j,q_1^{i+j-2s-2r}}^{-1} \bomega_{j,q_1^{
h-i-j+2s}},\end{equation*} where $h$ is the dual Coxeter number of
$\lie g$ if $\lie g$ is of type $B_n$ or $D_n$ and is twice the
dual Coxeter number if $\lie g$ is of type $C_n$. Given
$\boj\in\bj_r$, set
\begin{alignat*}{2}
\bpi_r(\boj) = &\ \bomega_{r,q_1^{r-i}}\prod_{s=1}^k \bpi_r(j_s-1,s-1)\bpi^{-1}_r(j_s,s),
&& {\text{if\ either \ }} \ \lie g= C_n {\text{ or }}  j_k<\bar n,\\
=&\ \bomega_{r,q_1^{r-i}}\left(\prod_{s=1}^{k-1}\bpi_r(j_s-1,s-1)\bpi^{-1}_r(j_s,s)\right)
\bpi_r(n-1,k-1) \times\\
&\ \times \ \bpi^{-1}_r(n,k-{1}/{4}), \
&& {\text{if}} \ \lie g= B_n {\text{ and }}\  j_k =n,\\
= &\ \bomega_{r,q^{r-i}}\left
(\prod_{s=1}^k
\bpi_r(j_s-1,s-1)\bpi^{-1}_r(j_s,s)\right)\bpi^{-1}_r(n,
k+{1}/{2}),&& {\text{if}} \ \lie g= D_n {\text{ and }}\  j_k =n-1.
\end{alignat*}
  where $\bar n= n$ (resp. $\bar n=n-1$) if $\lie g=B_n$
(resp. if $\lie g=D_n$). We understand that if $\boj$ is the empty
partition, then $\bpi_r(\boj)=\bomega_{r,q_1^{r-i}}$ and, if $r=0$, that $\bomega_{0,a}=1$ and $\omega_0=0$.

If $\lie g$ is of type $B_n$, define $\bpi_r(\boj,*)\in\cal P_q$
by
\begin{alignat*}{2}
\bpi_r(\boj,*)=&\ \bomega_{n,q^{2(n-i+2k)-1}}\bomega^{-1}_{n,
q^{2(n+i-2k-2r)-1}}, &&\text{ if } k \ne M  \text{ and } j_k\ne
n,\\
 =&\ 1,&& {\text{ otherwise} }.
\end{alignat*}

If $\lie g$ is of type $D_n$, define elements $\bpi_r(\boj,\pm)\in\cal P_q$ by
\begin{align*}
\bpi_r(\boj,+)&=\bomega_{n, q^{n-i+2k-1}} \bomega^{-1}_{n,
q^{n+i-2r-2k-1}},\\  \bpi_r(\boj,-)&=\bomega_{n-1, q^{n-i+2k-1}}
\bomega^{-1}_{n-1, q^{n+i-2r-2k-1}}.
\end{align*}

\vfill\eject

\begin{thm}\label{main}Let $V=V(\bomega_{i,1})$.

\begin{enumerit}
\item If $\lie g$ is of type $C_n$, the
assignment $\bj_r\to\cal P_q$ defined by $\boj\mapsto\bpi_r(\boj)$ is
injective and the image is $\wtl(V_{\omega_r})$. In particular,
$$\chl(V)=\sum_{r\in I_i}\sum_{w\in
W_{\omega_r}}\ \sum_{\boj\in\bj_r} T_w(\bpi_r(\boj)).$$
\item  If $\lie g$ is of type $B_n$, the
assignment $\bj_r\to\cal P_q$ defined by
$\boj\mapsto\bpi_r(\boj)\bpi_r(\boj,*)$ is injective and the image is
$\wtl(V_{\omega_r})$.  In particular, $$\chl(V)=\sum_{r\in I_i}\
\sum_{w\in W_{\omega_r}}\ \sum_{\boj\in\bj_r}
T_w(\bpi_r(\boj)\bpi_r(\boj,*)). $$
\item If $\lie g$ is of type $D_n$, then
$$\wtl(V_{\omega_r})=\{\bpi_r(\boj)\bpi_r(\boj,\pm):\boj\in\bj_{k,r}, 0\le k<M\} \cup \{\bpi_r(\boj):\boj\in\bj_{M,r}\}.$$
Moreover $$\chl(V)=\sum_{r\in I_i}\ \sum_{w\in
W_{\omega_r}}\ \left(\sum_{\boj\in\bj_r\backslash J_{M,r}}
\left(T_w(\bpi_r(\boj)\bpi_r(\boj,+))+
T_w(\bpi_r(\boj)\bpi_r(\boj,-))\right)
+\sum_{\boj\in\bj_{M,r}}T_w(\bpi_r(\boj))\right).$$

\end{enumerit}

\end{thm}

We  prove the theorem in the next three sections using Theorem
\ref{wconjfund} in an inductive way.

\section{The Case of $C_n$}

\subl{} Observe that the set $\bj_r$ depends on
$n,i,r$ and it will be necessary for the proofs to write $\bj_r$ as
$\bj_r(i)$. Notice moreover that
\begin{equation}\label{indeq}
(j_1,\cdots,j_M)\in\bj_r(i) \Leftrightarrow (j_2,\cdots ,j_M)\in
\bj_{j_1}(i-r+j_1-2),
\end{equation}
and also that
\begin{equation}\label{indeq2}  (j_1,\cdots ,j_{M-1})\in\bj_{r+2}(i)
\Leftrightarrow (j_1-2,j_2-2,\cdots j_{M-1}-2, j_M)\in \bj_r(i)\ \
\forall\ \ j_{M-1}-2<j_M<n.
\end{equation}
\begin{lem}\label{Cndim} We have
$$ |\bj_r(i)|  =\binom{n-r}{ M} - \binom{n-r}{M-1}= \dim
V(\omega_i)_{\omega_r}.$$
\end{lem}

\begin{pf} It suffices to prove the first equality, the second being well--known
 (see \cite{FH} for instance).
  If $M=0,1$ the result clearly holds for all
 $n\in\bn$ and $1\le i\le n$.
 Assume now that we know the result for all $M'<M$,  $n\in\bn$, and
  $1\le i\le n$.
By \eqref{indeq} we get,
$$|\bj_r(i)| = \sum_{j=r+1}^{n-2M+1} |\bj_{j_1}(i-r+j_1-2)| =
\sum_{j=r+1}^{n-2M+1} \left(
\binom{n-j}{M-1}-\binom{n-j}{M-2}\right),$$ where in the last
equality we used the induction hypothesis. The identity
$$\binom{m}{l} = \binom{m-1}{l}+\binom{m-1}{l-1}, \ \ 1\le l\le
m-1,$$ now gives the result.
\end{pf}

\subl{}

\begin{lem}\label{Cniota}
The map $\boj\mapsto\bpi_r(\boj)$ from $\bj_r(i)\to \cal P_q$ is
injective.\end{lem}
\begin{pf} We proceed by induction on $M$. Suppose that $\bpi_r(\boj)=\bpi_r(\boj')$ for
some $\boj,\boj'\in\bj_r(i)$. Writing $\boj=(j_1,\cdots ,j_M)$ and
$\boj'=(j_1',\cdots ,j_M')$ and comparing   the $(j_1-1)^{th}$
entries in $\bpi_r(\boj)$ and $\bpi_r(\boj')$ we find that
$j_1=j_1'$. This proves that induction starts at $M=1$. The
inductive step follows by \eqref{indeq}.
\end{pf}

\subl{} For $r>1$ and $r-1\le j< n$ define elements $w_{r,j}\in W$
and $T_{r,j}\in\cal B(\lie g)$  by \begin{align*} w_{r,j}& =
s_{j-1}s_{j-2}\cdots s_{r-1} s_{j+1}s_{j+2}\cdots s_{n-1} s_n\cdots s_r,\\
T_{r,j}&= T_{w_{r,j}}.
\end{align*} It is not hard to check (see \cite{H} for instance) that $w_{r,j}\in W_{\omega_r}$.

 The next proposition is   a straightforward if a
somewhat tedious  computation.
\begin{prop}\label{wrjCn} \hfill
\begin{enumerit}
\item For all $r\in I$, and $r-1\le j< n$, we have
$$w_{r,j}\omega_r = \omega_{r-2}+\alpha_j.$$
\item \begin{equation}\label{Cn0}
T_{r,j}(\bomega_{l,a}) =
\begin{cases}
\bomega_{l-2,aq^2} \bomega^{-1}_{j-1,aq^{j-l+3}}
\bomega_{j,aq^{j-l+2}} \bomega_{j,aq^{2n-j-l+2}}
\bomega^{-1}_{j+1,aq^{2n-j-l+3}}, & \text{ if } r\le l\le j,\\
\bomega_{l,aq^2} \bomega_{j,aq^{l-j}} \bomega_{j,aq^{2n-j-l+2}}
\bomega^{-1}_{j+1,aq^{l-j+1}} \bomega^{-1}_{j+1,aq^{2n-j-l+3}}, &
\text{ if } l> j.
\end{cases}
\end{equation}
\end{enumerit}\hfill\qedsymbol
\end{prop}

\subl{} Part (i) of Theorem \ref{main} now follows from Lemma
\ref{Cniota} and the next proposition.

\begin{prop}\label{Cnmain}
We have \begin{equation}\label{wtcn}\wtl(V_{\omega_r}) =
\{\bpi_r(\boj): \boj\in \bj_r\},\end{equation} and
\begin{equation}\label{multCn}
\dim V_{\bpi_r(\boj)} =1 \qquad \text{for all}\qquad \boj\in\bj_r.
\end{equation}
\end{prop}

\begin{pf}
Recall from\cite{CPbanff} that $$V \cong V(\omega_i)$$ as
$U_q(\lie g)$-modules. It suffices to prove that,
\begin{equation}\label{wtcn'}
\{\bpi_r(\boj): \boj\in \bj_r\} \subset \wtl(V_{\omega_r}),
\end{equation}
since Lemma \ref{Cndim} and Lemma \ref{Cniota}  then imply both \eqref{wtcn} and \eqref{multCn}.

To prove \eqref{wtcn'} we proceed by induction on $M$ with
induction beginning at $M=0$. Assume that we know the result for
$M-1$. To prove the inductive step it follows from \eqref{indeq2}
that if $\boj=(\boj_1,\cdots ,\boj_M)\in\bj_r$ then
$\boj'=(\boj_1+2,\cdots ,\boj_{M-1}+2)\in\bj_{r+2}$.
 The induction hypothesis implies
that $\bpi_{r+2}(\boj')\in \wtl(V)$ and hence by Theorem
\ref{wconjfund} we see that $T_{r+2,j}(\bpi_{r+2}(\boj')) \in \wtl
(V)$ for all $r<j<n$.

Using  Proposition \ref{wrjCn}(ii) we find that
  $$T_{r+2,j}(\bpi_{r+2}(\boj')) =
\bpi_{r}(\boj)\balpha_{j,q^{2n+i-j-2r-2M}},$$ and also that the
$j^{th}$--coordinate of $T_{r+2,j}(\bpi_{r+2}(\boj'))$ is
$$(1-q^{i+j-2r-2M}u)(1-q^{2n+i-j-2r-2M}u).$$ Hence by Theorem
\ref{wconjfund}(ii) we see that
$$\bpi_r(\boj)=T_{r+2,j}(\bpi_{r+2}(\boj'))\balpha^{-1}_{j,q^{2n+i-j-2r-2M}}\in\wtl(V_{\omega_r}).$$
\end{pf}

\section{The Case of $B_n$.}

\subl{}
\begin{lem}\label{cardBn} We have $$ |\bj_r|=
\sum\limits_{k=0}^{M}\binom{n-r}{k} = \dim V_{\omega_r}.$$\hfill\qedsymbol
\end{lem}

\begin{pf}
The first equality is clear. For the second, recall that it was
proved  in \cite{CPbanff} that $V \cong \opl_{l=0}^{[i/2]}
V(\omega_{i-2l})$
 as $U_q(\lie g)$-modules. Since  $$\dim V(\omega_j)_{\omega_r}
 = \binom{n-r}{[\frac{j-r}{2}]},\ \ j<n,$$ it follows that
 $$\dim V_{\omega_r} =
 \sum\limits_{l=0}^{M}\binom{n-r}{[\frac{i-2l-r}{2}]} =
 \sum\limits_{k=0}^{M}\binom{n-r}{k}.$$
\end{pf}

\subl{}
\begin{lem}\label{Bn1-1}
The map $\boj \mapsto \bpi_r(\boj)\bpi_r(\boj,*)$ from $\bj_r$ to
$\cal P_q$ is injective.
\end{lem}

\begin{pf}
Suppose that $0\le k,k'\le M$,  $\boj,\boj'\in\bj_r$,
$\boj=(j_1,\cdots, j_k)$, $\boj'=(j'_1,\cdots, j'_{k'})$ are such
that
\begin{equation}\label{Bninjn}
\bpi_r(\boj)\bpi_r(\boj,*)=\bpi_r(\boj')\bpi_r(\boj',*).
\end{equation}
We first show that  $k=k'$. For this, notice that for any
$\boj''=(j_1'',\cdots ,j_{k''}'')\in\bj_r$ we have
$$(\bpi_r(\boj'')\bpi_r(\boj'',*))_n =1\iff k''=M \ \ {\text{and}}\ \
j''_{k''}<n.$$ Hence,  $(\bpi_r(\boj)\bpi_r(\boj,*))_n =1$ implies
$k=k'=M$. Otherwise, the equation
$$(\bpi_r(\boj)\bpi_r(\boj,*))_n=(\bpi_r(\boj')\bpi_r(\boj',*))_n$$
gives that either   $j_k, j'_{k'}<n$ or   $j_k=j'_{k'}=n$. In
the first case we get $\bpi_r(\boj,*)_n= \bpi_r(\boj',*)_n$ and in
the second case  we get $ \bpi_r(\boj)_n= \bpi_r(\boj')_n$. In any
case it follows that  $k=k'$.

 Now, suppose $\boj\ne \boj'$ and let $1\le s_0\le n$ be
minimal such that $j_{s_0}\ne j_{s_0}'$.  Assume without loss of generality that $
j'_{s_0}>j_{s_0}$. This means that $
j'_{s_0-1}=j_{s_0-1}<j'_{s_0}-1$ and hence
\begin{equation}\label{js0}\bpi_r(\boj')_{j'_{s_0}-1}
=
(1-q_1^{2n-i-j'_{s_0}+2s_0-2}u)(1-q_1^{i+j'_{s_0}-2s_0-2r+1}u)^{-1}
\ne 1.\end{equation} We claim now that there exists $s_1\ge s_0$
such that $j_{s_1}=j'_{s_0}-1$. Assuming the claim, we get a
contradiction to the fact that $\boj\ne\boj'$ as follows. Since
$$\bpi_r(\boj')_{j'_{s_0}-1}
=(1-q_1^{i+j'_{s_0}-2s_1-2r-1}u)(1-q_1^{2n-i-j'_{s_0}+2s_1}u)^{-1}=\bpi_r(\boj)_{j_{s_1}},$$
comparing with \eqref{js0} gives $$2n-i-j'_{s_0}+2s_0-2 =
i+j'_{s_0}-2s_1-2r-1,$$ which is obviously  impossible.
 To  prove the claim, set $$s_1=\max\{1\le s\le
k:j_s<j'_{s_0}\}.$$ The claim follows if we prove that $j_{s_1+1}>
j'_{s_0}$. The maximality of
$s_1$ implies that $j_{s_1+1}\ge j'_{s_0}$ and hence it suffices
to prove that $j_{s_1+1}\ne j'_{s_0}$. If $j_{s_1+1}=j'_{s_0}$
then the same argument that gave \eqref{js0} gives,
$$\bpi_r(\boj)_{j_{s_1+1}-1}=(1-q_1^{2n-i-j_{s_1+1}+2(s_1+1)-2}u)(1-q_1^{i+j'_{s_1+1}-2(s_1+1)
_-2r+1}u)^{-1} \ne 1 ,$$ which implies that $s_1+1=s_0$
contradicting $s_1\ge s_0$.
\end{pf}

\subl{}

 For $r>0$ and $r\le j<n$ define elements $w_{r,j}\in
 W_{\omega_r}$ by,
\begin{align*}w_{r+1,j} &= s_{j-1}s_{j-2}\cdots s_{r} s_{j+1}s_{j+2}\cdots
s_{n-1} s_n\cdots s_{r+1},\\
w_{r,n}&=s_{n-1}\cdots s_{r},\\  T_{r,j} &=
T_{w_{r,j}}.\end{align*} The proof of the next proposition is
along the same lines as the proof of Proposition \ref{wrjCn} and we
omit the details.

\begin{prop}\label{wrjBn}\hfill
\begin{enumerit}
\item[(i)] For all $r\in I\backslash \{n\}$ and $r-1\le j< n$ we have
$$w_{r,j}\omega_{r} = \omega_{r-2}+\alpha_j \quad\text{ and }\quad w_{r,n}\omega_r=\omega_{r-1}+\alpha_n.$$

\item[(ii)] For $r-1\le j<n-1$ and $r\le l<n$ we have:
\begin{align}\label{Bnjl}
T_{r,j}(\bomega_{l,1}) =
\begin{cases}
\bomega_{l-2,q_1^2} \bomega^{-1}_{j-1,q_1^{j-l+3}}
\bomega_{j,q_1^{j-l+2}} \bomega_{j,q_1^{2n-j-l-1}}
\bomega^{-1}_{j+1,q_1^{2n-j-l}},
 &\text{ if }\ l\le j,\\
\bomega_{l,q_1^2} \bomega_{j,q_1^{l-j}} \bomega_{j,q_1^{2n-j-l-1}}
\bomega^{-1}_{j+1,q_1^{l-j+1}} \bomega^{-1}_{j+1,q_1^{2n-j-l}},
&\text{ if }\  l> j.
\end{cases}
\end{align}
Further,
\begin{align}\notag
T_{r,j}(\bomega_{n,q^{-1}})& = \bomega_{j,q_1^{n-j-1}}
\bomega^{-1}_{j+1,q_1^{n-j}} \bomega_{n,q^3},\\ \label{Bnnl}
T_{r,n-1}(\bomega_{l,1})&= \ \bomega_{l-2,q_1^2}
\bomega^{-1}_{n-2,q_1^{n-l+2}} \bomega_{n-1,q_1^{n-l+1}}
\bomega_{n-1,q_1^{n-l}} \bomega^{-1}_{n,q^{2(n-l)+1}}
\bomega^{-1}_{n,q^{2(n-l)+3}},\\ \notag
T_{r,n-1}(\bomega_{n,q^{-1}})& = \ \bomega_{n-1,1}
\bomega^{-1}_{n,q},\\ \notag T_{r,n} (\bomega_{l,1}) = &\
\bomega_{l-1,q_1}\bomega^{-1}_{n-1,q_1^{n-l+1}}\bomega_{n,q^{2(n-l)-1}}\bomega_{n,q^{2(n-l)+1}}.
\end{align}

\end{enumerit}\hfill\qedsymbol
\end{prop}

\subl{} To prove Theorem \ref{main} (ii)
 we proceed  by induction on $M$. Induction clearly begins when $M=0$.
 The inductive step
is immediate from the following proposition, Lemma
\ref{cardBn}, and Lemma \ref{Bn1-1}.

\begin{prop} Assume that $M>0$ and let $\boj=(j_1,\cdots, j_k)\in\bj_{r}$.
\begin{enumerit}
\item If $k<M$ and $j_k<n$ we have:
\begin{equation*}
\bpi_r(\boj)\bpi_r(\boj,*) = T_{r+1,n}(\bpi_{r+1}(\boj')\bpi_{r+1}(\boj',*))\balpha^{-1}_{n,a}
\in\wtl(V_{\omega_r}),
\end{equation*}
where  $\boj' = (j_1+1,\cdots, j_k+1)\in \bj_{r+1}$ and
$a=q^{2(n+i-2r-2k)-3}$.
\item If $j_k=n$, we have \begin{equation*}
\bpi_r(\boj)\bpi_{r}(\boj,*) = T_{r+1,n}(\bpi_{r+1}(\boj')\bpi_{r+1}(\boj',*))\balpha^{-1}_{n,a}
\in\wtl(V_{\omega_r}),
\end{equation*}
where  $\boj'=(j_1+1,\cdots, j_{k-1}+1)\in\bj_{r+1}$ and
$a=q^{2(n-i+2k)-5}$.\item  If $k=M$ and $j_k <n$, then
\begin{equation*}
\bpi_r(\boj)\bpi_{r}(\boj,*) = T_{r+2,j}(\bpi_{r+2}(\boj')\bpi_{r+2}(\boj',*))\balpha^{-1}_{j,a}
\in\wtl(V_{\omega_r}),
\end{equation*}
where $\boj'=(j_1+2,\cdots, j_{k-1}+2)\in\bj_{r+2}$, $
a=q_1^{2n-j-r-3}$, and $j= j_k$.
\end{enumerit}
\end{prop}

\begin{pf} Observe first that it is clear that the elements
$\boj'$ defined in the proposition are in $\bj_{r+1}$ in the first
two cases and in $\bj_{r+2}$ in the third case. The fact that
$\bpi_r(\boj)\bpi_{r}(\boj,*)$ and $\bpi_{r+1}(\boj')\bpi_{r+1}(\boj',*)$ (resp.
$\bpi'_{r+2}(\boj')\bpi_{r+2}(\boj',*))$ are related as in the proposition is again a
tedious but simple checking using the formulas in Proposition
\ref{wrjBn}(ii). The main point is to notice that this implies
 $\bpi_r(\boj)\bpi_{r}(\boj,*)\in\wtl(V_{\omega_r})$. For that, one observes
that the calculation gives respectively:

\begin{align*}\hspace*{-20pt}{\rm(i)} \ \ (T_{r+1,n}(\bpi_{r+1}(\boj')&\bpi_{r+1}(\boj',*)))_n = \\
= &\ (1-q^{2(n+i-2r-2k)-3}u)(1-q^{2(n-i+2k)-1}u),\ \ {\text{if}}\
j_k=n-1 \text{ or }\ \ k<M-1/2,\\
=& \ (1-q^{2(n+i-2r-2k)-3}u)(1-q^{2(n+i-2r-2k)-5}u) {\text{ if}}\
j_k<n-1 \text{ and } \ \ k=M-1/2.
\end{align*}

\item[(ii)] $(T_{r+1,n}(\bpi_{r+1}(\boj')\bpi_{r+1}(\boj',*)))_n =  (1-q^{2(n-i+2k)-5}u)(1-q^{2(n+i-2r-2k)+1}u)$.
\vskip 6pt

\item[(iii)] $(T_{r+2,j}(\bpi_{r+2}(\boj')\bpi_{r+2}(\boj',*)))_j =  (1-q_1^{j-r}u)(1-q_1^{2n-j-r-3}u)$.

\vskip15pt
The result then follows from Theorem \ref{wconjfund}.
\end{pf}

\section{The Case of $D_n$}

\subl{}
\begin{lem} We have:
$$|\bj_{M,r}|+2\sum_{k=0}^{M-1}|\bj_{k,r}| =
\sum\limits_{l=0}^{M} \binom{n-r}{l} = \dim V_{\omega_r}.$$
\end{lem}

\begin{pf}
It follows from  the definition of $\bj_{k,r}$ that $$|\bj_{k,r}|
= \binom{n-r-1}{k},$$ and hence to prove the  the first equality
we must show that $$\sum\limits_{l=0}^{M} \binom{n-r}{l} =
\binom{n-r-1}{M} + 2\sum\limits_{k=0}^{M-1} \binom{n-r-1}{k}.$$
Using the binomial identity $$\binom{n-r}{l} = \binom{n-r-1}{l} +
\binom{n-r-1}{l-1},$$ we find that
\begin{align*}
\sum\limits_{l=0}^{M} \binom{n-r}{l} &=\ 1+\sum\limits_{l=1}^{M}
\left(\binom{n-r-1}{l}+\binom{n-r-1}{l-1}\right) \\ &=
1+\binom{n-r-1}{M} + \sum\limits_{l=1}^{M-1}
\binom{n-r-1}{l}+\sum\limits_{l=0}^{M-1} \binom{n-r-1}{l}\\ &=\
\binom{n-r-1}{M} + 2\sum\limits_{l=0}^{M-1} \binom{n-r-1}{l}.
\end{align*}

For the second equality, recall that it was proved in \cite{CPbanff}
that as $U_q(\lie g)$-modules $$V \cong \opl_{l=0}^{[i/2]}
V(\omega_{i-2l}).$$ The result now follows since $$\dim
V(\omega_j)_{\omega_r} = \binom{n-r}{(j-r)/2},\ \ 1\le j\le n-2.$$
\end{pf}

\subl{} Given $r>1$ and $r-1\le j\le n$, define elements
$w_{r,j}\in W_{\omega_r}$ by,
\begin{align*} w_{r,j}&=s_{j-1}s_{j-2}\cdots s_{r-1}s_{j+1}\cdots
s_{n-2}s_ns_{n-1}\cdots s_r,\ \ j\le n-2,\\ &=s_{n-2}\cdots
s_{r-1}s_{j'} s_{n-2}\cdots s_r,\ \ j,j'\in\{n-1, n\},\ \ j'\ne
j,\\ T_{r,j}&=T_{w_{r,j}}.\end{align*}

\begin{prop}\label{wrjDn}For all $1< r\le n-2$ and $r-1\le j\le n$ we have:\hfill
\begin{enumerit}
\item[(i)] $w_{r,j}\omega_{r} = \omega_{r-2}+\alpha_j$.
\item[(ii)]
\begin{align}
T_{r,j}(\bomega_{l,a}) =
\begin{cases}
\bomega_{l,aq^2} \bomega_{j,aq^{l-j}}\bomega_{j,aq^{2n-l-2-j}}\bomega^{-1}_{j+1,aq^{l-j+1}}\bomega^{-1}_{j+1,aq^{2n-l-j-1}}, & \text{if } j<l,\\
\bomega_{l-2,aq^2}\bomega_{j-1,aq^{j-l+3}}^{-1}\bomega_{j,aq^{j-l+2}}\bomega_{j,aq^{2n-l-j-2}}\bomega_{j+1,aq^{2n-l-j-1}}^{-1}, & \text{if } l\leq j\leq n-2,\\
\bomega_{l-2,aq^2}\bomega_{n-2,aq^{n-l+2}}^{-1}\bomega_{j,aq^{n-l-1}}\bomega_{j,aq^{n-l+1}}, & \text{if } j= n-1,n,
\end{cases}
\end{align}
if $1\le l\le n-2$ and
\begin{align}
T_{r,j}(\bomega_{l,a}) =
\begin{cases}
\bomega_{j,aq^{n-1-j}}\bomega_{j+1,aq^{n-j}}^{-1}\bomega_{l',aq^2}, & \text{if } j<n-2,\\
\bomega_{n-2,aq}\bomega^{-1}_{l,aq^{2}}, & \text{if }  j= n-2,\\
\bomega_{l,a}, & \text{if } j= l,\\
\bomega_{n-3,aq^{2}}\bomega_{n-2,aq^{3}}^{-1}\bomega_{l',aq^2}, & \text{if } j=l',
\end{cases}
\end{align}
if $l=n-1,n$, where $l'\in\{n-1,n\}\backslash\{l\}$.
\end{enumerit}\hfill\qedsymbol
\end{prop}

\subl{} The next proposition is proved in a similar manner to the
corresponding one for  $B_n$ and $C_n$. We omit the details this
time.

\begin{prop}
 For $M\ge 0$, we have
$$ \{\bpi_r(\boj):\boj\in\bj_{M,r}\} \cup \{\bpi_r(\boj)\bpi_r(\boj,\pm):\boj\in\bj_{k,r}, 0\le k<M\}\subset
\wtl(V_{\omega_r}).$$\hfill\qedsymbol\end{prop} \subl{} To complete the proof of Theorem
\ref{main}(iii), we must prove that in fact
$$ \{\bpi_r(\boj):\boj\in\bj_{M,r}\} \cup \{\bpi_r(\boj)\bpi_r(\boj,\pm):\boj\in\bj_{k,r}, 0\le k<M\}=
\wtl(V_{\omega_r}).$$For
$D_n$ this is more difficult, since it is no longer true that the
maps $\boj\mapsto \bpi_r(\boj)\bpi_r(\boj,\pm)$ are injective.

The next lemma is a simple checking.
\begin{lem}\label{dist} Let $\boj\in\bj_{k,r}$ and $\boj'\in \bj_{k',r}$.
\begin{enumerit}
\item  We have $\bpi_r(\boj, \pm) = \bpi_r(\boj',\pm)$ iff $k=k'$. Moreover, if $k=M$, then $\bpi_r(\boj, \pm) =1.$
\item If $\bpi_r(\boj)\bpi_r(\boj,+) = \bpi_r(\boj')\bpi_r(\boj',\pm)$, then  $k=k'$. Moreover, if  $k=k'<M$, then
$\bpi_r(\boj)\bpi_r(\boj,+) \ne \bpi_r(\boj')\bpi_r(\boj',-)$.
\end{enumerit} \hfill\qedsymbol
\end{lem}

\subl{} Define an equivalence relation $\sim$ on $\bj_{k,r}$ by
$$ \boj\sim \boj' \iff \ \bpi_r(\boj)=\bpi_r(\boj').$$
Let $\bar \boj$ be the equivalence
class of $\boj$.

\begin{prop}\label{complete} Let $\boj\in\bj_{k,r}$. Then $$\dim
V_{\bpi_r(\boj)\bpi_r(\boj,\pm)}\ge |\bar\boj|.$$\end{prop}

Assuming this proposition the proof of Theorem \ref{main}(iii) is
completed as follows. For each $0\le k\le M$, fix a set
$\bs_{k,r}\subset\bj_{k ,r}$ of representatives of the distinct
equivalence classes with respect to $\sim$.  Proposition
\ref{complete} and Lemma \ref{dist} imply that
\begin{align*}
\sum_{\boj\in\bs_{M,r}} \dim(V_{\bpi_r(\boj)}) & + \
\sum_{k=0}^{M-1} \sum_{\boj\in\bs_{k,r}} \left(\dim(V_{\bpi_r(\boj)\bpi_r(\boj,+)} ) +  \dim(V_{\bpi_r(\boj)\bpi_r(\boj,-)} ) \right) \ge\\
 & \ \sum_{\boj\in\bs_{M,r}}|\bar\boj| + 2\sum_{k=0}^{M-1}\sum_{\boj\in\bs_{k,r}}|\bar\boj|  = |\bj_{M,r}|+2\sum_{k=0}^{M-1}|\bj_{k,r}|.
 \end{align*}

Since $V_{\bpi_r(\boj)\bpi_r(\boj,\pm)}\subset V_{\omega_r}$ it
follows from Lemma 5.1 that $$V_{\omega_r}=
\left(\bigoplus_{\boj\in\bs_{M,r}}  (V_{\bpi_r(\boj)})\right)
\oplus \left(\bigoplus_{k=0}^{M-1} \bigoplus_{\boj\in\bs_{k,r}}
(V_{\bpi_r(\boj)\bpi_r(\boj,+)} \oplus
V_{\bpi_r(\boj)\bpi_r(\boj,-)})\right),$$ which proves Theorem
\ref{main}(iii).

\subl{} It remains to prove Proposition \ref{complete}. This
requires some combinatorial definitions and results which we now
establish and which allow us to actually compute $|\bar\boj|$ in
terms of $\boj\in\bj_{k,r}$. Thus we shall see that
\begin{equation}|\bar\boj|=2^{M_\boj},\end{equation} where $M_\boj=|\{s:1\le s\le k,
j_s\in\{n-i+r+2s-2, n-i+r+2s-1\}\}|.$

\subl{} Recall that a strictly increasing partition $\bon$  of
positive integers  of length $k$ is an increasing sequence  $0<
j_1<j_2<\cdots< j_k$ of natural numbers. Let
$$\supp(\bon)=\{j_s:1\le s\le k\},$$ and let
$\iota_\bon:\supp(\bon)\to \bn$ be defined by $\iota_\bon(j)=s$,
if $j=j_s$ for $ 1\le s\le k$. Clearly any finite subset of $\bn$
defines a strictly increasing partition.

Given a finite subset $S\subset\bn$ and  a partition $\bon$ of
length $k$, let $\bon_S$ be the partition corresponding to the
set, $$(\supp(\bon)\ \backslash S) \cup (S\ \backslash
\supp(\bon)).$$ Clearly, $$ \bon_S\ne \bon_{S'} \ \ \text{if}\ \
S\ne S',$$ and
\begin{equation}\label{invpart}
(\bon_S)_{S}=\bon, \ \ (\bon_S)_{S'}=(\bon_{S'})_S=(\bon)_{S\cup
S'}, \ \ \text{if } S\cap S'=\emptyset.
\end{equation}

We shall adopt the convention that if $k,k'\in\bz$, then
$[k,k']=\{\min\{k,k'\},\min\{k,k'\}+1,\cdots ,\max\{k,k'\}\}$. Define also $(k,k']$ and $[k',k)$ in the obvious way.

Associated with a partition
$\bon$, define functions $\sigma_\bon^\pm,
\tau^\pm_{\bon}:\supp(\bon)\to \bn$ by,
\begin{align*}
\sigma^+_\bon(j)& = \max\{j'\in\supp(\boj): j'\ge j \text{ and } j''-j< 2(\iota_\bon(j'')-\iota_\bon(j)) \ \forall\ j''\in\supp(\boj)\cap (j,j']\},\\
\sigma^-_\bon(j)& = \min\{j'\in\supp(\boj): j'\le j \text{ and } j-j''< 2(\iota_\bon(j)-\iota_\bon(j'')) \ \forall\ j''\in\supp(\boj)\cap [j',j)\},
\end{align*}
and
$$\tau^\pm_{\bon}(j)=j+2(\iota_\bon(\sigma_\bon^{\pm}(j))-\iota_\bon(j))\pm 1.$$

 The
following lemma is easy.
\begin{lem}\label{cardinvpart} Let $j\in\supp(\bon)$. Then,
$$|\supp(\bon_{[j,\tau^{\pm}_{\bon}(j)]})| =
|\supp(\bon)|.$$\hfill\qedsymbol
\end{lem}

>From now on we set \begin{equation}
\bon^\pm(j)=\bon_{[j,\tau^\pm_{\bon}(j)]},\ \
j\in\supp(\bon).\end{equation} Moreover, if $S\subset\supp(\bon)$
is such that \begin{equation}\label{disjoint}
[j,\tau^\pm_{\bon}(j)]\cap [j',\tau^\pm_{\bon}(j')]=\emptyset,\ \
\forall j\ne j', j,j'\in S,\end{equation} then set
\begin{equation} \bon^\pm(S) =(\bon^\pm(j))^\pm(S\backslash
\{j\}),\ \ j\in S\subset \supp(\bon).
\end{equation}
Notice that $\bon^\pm(S)$ is well-defined by \eqref{invpart}.

\subl{} For a partition $\bon$ and an integer $m> 0$, define
\begin{align*} \supp_m(\bon)&
= \{j\in\supp(\bon): 2\iota_\bon(j)=j-m\}.
\end{align*}

The following lemma is easy.

\begin{lem}\label{cardcomp} Let $j\in\supp_m(\bon)$ for some $m>1$.
Then $\tau^\pm_\bon(j)\in \supp_{m\pm 1}(\bon^{\pm}(j))$, and
$$\supp_{m\pm 1}(\bon^{\pm}(j)) = \supp_{m\pm 1}(\bon)\sqcup
\{\tau^\pm_\bon(j)\},\ \ \ \ \ \supp_m(\bon^{\pm}(j)) =
\supp_m(\bon)\backslash \{j\}.$$ In particular
$\tau^{\mp}_{\bon^{\pm}(j)}(\tau^\pm_\bon(j))=j$.\hfill\qedsymbol
\end{lem}
\begin{rem} Notice that $\bon^{\pm}(S)$ is well defined for all $S\subset \supp_m(\boj)$. The Lemma then implies that
$\supp_m(\bon^\pm(\supp_m(\bon))=\emptyset$.\end{rem}

\subl{} From now on we set
\begin{equation}N = n-i+r-2,\ \ \ \supp^+(\boj)=\supp_N(\boj),\ \
\ \supp^-(\boj)=\supp_{N+1}(\boj).
\end{equation} Let $\boj\in \bj_{k,r}$.
 The next
proposition describes $\bar \boj$ for $\boj\in \bj_{k,r}$.
\begin{prop}\label{Dncomb}
Let $\boj,\boj'\in\bj_{k,r}$. Then, \begin{equation}
\label{eq}\boj\sim\boj'\ \iff
\boj'=(\boj^-(S^-))^+(S^+),\end{equation} for some
  $S^-\subset
\supp^-(\boj)$ and  $S^+\subset \supp^+(\boj^-(S^-))$. In
particular, \begin{equation}\label{cardboj} |\bar\boj|=
2^{|\supp^+(\boj)|+|\supp^{-}(\boj)|}.\end{equation}
\end{prop}

\begin{cor} Let $\boj\in\bj_{k,r}$ be such that
$\supp^-(\boj)=\emptyset$. Then $|\boj|=2^{|\supp^+(\boj)|}$ and
$$|\bj_{k,r}|=\sum_{\substack{\boj\in\bj_{k,r},\\
\supp^-\boj=\emptyset}}2^{|\supp^+(\boj)|}.$$
\hfill\qedsymbol\end{cor}

Notice that \eqref{cardboj} is immediate from \eqref{eq} and Lemma
\ref{cardcomp}.

For the proof, it is useful to notice that,
\begin{equation}\label{vanish}
\bpi_r(j,s) = 1 \Leftrightarrow 2s=j-(N+1).
\end{equation}

\subl{} We now see that  Proposition \ref{Dncomb} can be deduced
from the following.
\begin{prop}\label{sim}
Let $\boj\in\bj_{k,r}$. For any $S\subset\supp^\pm(\boj)$ we have $\boj^{\pm}(S) \sim \boj$.
\end{prop}
 If $\boj',\boj\in\bj_{k,r}$ are related as in the right hand side of \eqref{eq}, we
have by Proposition \ref{sim} that, $$\boj\sim\boj^-(S^-)
\sim(\boj^-(S^-))^+(S^+) = \boj'.$$ For the converse, let
$\boj\sim\boj'$ and assume that
\begin{equation}\label{uniquemin}
\boj^{-}(\supp^-(\boj)) = \boj'^{-}(\supp^-(\boj').
\end{equation}
Set $$S^- = \supp^-(\boj) \qquad\text{and}\qquad S^+ =
\tau^-_{\boj'}(\supp^-(\boj')).$$ It follows from Lemma
\ref{cardcomp} that $$S^+\subset
\supp^+(\boj'^{-}(\supp^-(\boj'))).$$ Then, using \eqref{invpart},
Lemma \ref{cardcomp} and \eqref{uniquemin} we see that $$\boj'=
\left(\boj'^{-}(\supp^-(\boj')) \right)^+(S^+) =
(\boj^{-}(S^-))^+(S^+).$$  It remains to show that
\eqref{uniquemin} is always satisfied if $\boj\sim\boj'$. Observe
that  Proposition \ref{sim} gives $$\boj^-(\supp^-(\boj)) \sim
\boj'^{-}(\supp^-(\boj')),$$ and that by
 Remark \ref{cardcomp} $$\supp^-(\boj^-(\supp^-(\boj)))
= \supp^-(\boj'^{-}(\supp^-(\boj')))= \emptyset.$$ In other words,
to prove \eqref{uniquemin} it suffices to prove,\begin{equation}
\boj_1\sim\boj_2, \ \text{and}\ \
\supp^-(\boj_1)=\supp^-(\boj_2)=\emptyset\implies
\boj_1=\boj_2.\end{equation} Indeed, if  $\boj_1\ne \boj_2$  set
$$j=\max\{j': j'\in \supp(\boj_1)\cup \supp(\boj_2) \text{ and }
j'\notin \supp(\boj_1)\cap \supp(\boj_2)\}.$$ Assume
$j\in\supp(\boj_1)$. Then, either $j+1\in \supp(\boj_1)\cap
\supp(\boj_2)$ or $j+1\notin \supp(\boj_1)\cup \supp(\boj_2)$. In
any case it follows that $$\bpi_r(\boj_1)_j = \bpi_r(\boj_2)_j =
1.$$ If $j+1\notin \supp(\boj_1)\cup \supp(\boj_2)$ then
$\bpi_r(\boj_1)_j = \bpi_r(j,\iota_{\boj_1}(j))$ and  \eqref{vanish} gives that $j\in \supp^-(\boj_1)$. Otherwise we have $\bpi_r(\boj_2)_j =
\bpi_r(j,\iota_{\boj_1}(j+1)-1)$ and
\eqref{vanish}  shows that $j\in \supp^-(\boj_1)$ in this case as well,
contradicting $\supp^-(\boj_1) =\emptyset$. This completes the
proof of Proposition \ref{Dncomb}.

\subl{}{\it Proof of Proposition \ref{sim}}.

It clearly suffices to prove the result when $S=\{j\}$ since the
general case follows by transitivity. We can assume that $j\in
\supp^+(\boj)$ since then using Lemma \ref{cardcomp}  the case
$j\subset\supp^-(\boj)$ follows. Since, $$\supp(\boj)\cap [r,j-1]
= \supp(\boj^+(S))\cap [r,j-1] \quad\text{and}\quad
\supp(\boj)\cap [\tau_\boj^+(j)+1,n] = \supp(\boj^+(S))\cap
[\tau_\boj^+(j)+1,n],$$ it follows immediately  that
$$\bpi_r(\boj)_{l} = \bpi_r(\boj^+(S))_{l} \ \forall\ l\in
[r,j-1)\cup  (\tau_\boj^+(j),n].$$ If $l\in [j-1,\tau^+_\boj(j)]$
we proceed by induction on $l$. To see that induction
starts at $l=j-1$, observe that $$j\notin \supp(\boj^+(S)),\ \
2(\iota_\boj(j)-1)=(j-1) - (N+1).$$ The conclusion follows from
\eqref{vanish}. For the inductive step, assume that
$$\bpi_r(\boj)_{l' } = \bpi_r(\boj^+(j))_{l'},\ \forall \ j-1\le
l'<l \le \tau_\boj^+(j).$$  Assume first that  $l\in \supp(\boj)$,
in particular $l<\tau_\boj^+(j)$. If $(l+1)\in \supp(\boj)$ then,
 $$\bpi_r(\boj)_l=\bpi_r(\boj^+(S))_l=1$$ and we are done. If
 $(l+1)\notin\supp(\boj)$, then
$(l+1)\in\supp(\boj^+(j))$ and we have from the definition of
$\bpi_r$ that, $$\bpi_r(\boj)_l = \bpi_r(l,\iota_\boj(l))_l
\qquad\text{ and }\ \ \bpi_r(\boj^+(j))_l =
\bpi_r(l,\iota_{\boj^+(j)}(l+1)-1)_l,$$ which gives,
$$\bpi_r(\boj)_l
=
(1-q^{i-2r+l-2\iota_\boj(l)}u)(1-q^{2n-i-l+2\iota_\boj(l)-2}u)^{-1},$$
and $$\bpi_r(\boj^+(j))_l =
(1-q^{2n-i-l+2\iota_{\boj^+(j)}(l+1)-4}u)(1-q^{i+l-2r-2\iota_{\boj^+(j)}(l+1)+2}u)^{-1}.$$
Let $$l_0 = \min\{m\ge j: m\le l'' \le l \implies l'' \in
\supp(\boj)\}.$$  Then, either   $l_0=j$, in which case we have
$$\iota_\boj(l) = \iota_\boj(j) +l-j = (i-r-n-j+2l+2)/2,\ \ \ \
\iota_{\boj^+(j)}(l+1)=\iota_\boj(j), $$  which implies $$
\bpi_r(\boj)_{l} = \bpi_r(\boj^+(j))_{l},$$ and we are done, or
$l_0>j$. In that case $(l_0-1)\in \supp(\boj^+(j))$,
$\iota_\boj(l) = \iota_\boj(l_0) + l-l_0,$ and
$$\bpi_r(\boj)_{l_0-1} = \bpi_r(l_0-1,\iota_\boj(l_0)-1)_{l_0-1},
\qquad\text{ and }\qquad \bpi_r(\boj^+(j))_{l_0-1} =
\bpi_r(l_0-1,\iota_{\boj^+(j)}(l_0-1))_{l_0-1}.$$ Since $l_0-1<l$,
the induction hypothesis gives $\bpi_r(\boj)_{l_0-1} =
\bpi_r(\boj^+(j))_{l_0-1}$, which implies that
$$\iota_{\boj^+(j)}(l_0-1) = i+l_0-n-r-\iota_\boj(l_0)+1$$ and,
therefore, $$\iota_{\boj^+(j)}(l+1) =
i+l_0-n-r-\iota_\boj(l_0)+2.$$ The conclusion follows.

Finally, suppose that
 $l\notin\supp(\boj)$. If $l=\tau_\boj^+(j)$, then \eqref{vanish} and Lemma \ref{cardinvpart}
 imply $$\bpi_r(\boj)_l = \bpi_r(\boj^+(j))_l =1,$$ and we are done.
 If $l<\tau_\boj^+(j)$ and $l+1\in\supp(\boj^+(j))$
  we again have $$\bpi_r(\boj)_l = \bpi_r(\boj^+(j))_l =1$$ and we are done.
   Otherwise we have $(l+1)\in\supp(\boj)$, and
$$\bpi_r(\boj)_l = \bpi_r(l,\iota_\boj(l+1)-1)_l ,\ \qquad
\bpi_r(\boj^+(j))_l = \bpi_r(l,\iota_{\boj^+(j)}(l))_l,$$ which
gives, $$\bpi_r(\boj)_l =
(1-q^{2n-i-l+2\iota_\boj(l+1)-4}u)(1-q^{i-2r+l-2\iota_\boj(l+1)
+2}u)^{-1}$$ and $$\bpi_r(\boj^+(j))_l =
(1-q^{i+l-2r-2\iota_{\boj^+(j)}(l)}u)(1-q^{2n-i-l+2\iota_{\boj^+(j)}(l)-2}u).$$
Let $$l_0 = \max\{l': j\le l'< l,  l'\in \supp(\boj) \}.$$ Then,
$l_0+1\in \supp(\boj^+(j))$ and  $$\iota_{\boj^+(j)}(l) =
\iota_{\boj^+(j)}(l_0+1) + l-(l_0+1), \qquad\qquad
\iota_\boj(l_0)= \iota_\boj(l+1)-1.$$ In particular,
$$\bpi_r(\boj)_{l_0} = \bpi_r(l_0,\iota_\boj(l_0))_{l_0} \qquad\text{and}
\qquad \bpi_r(\boj^+(j))_{l_0} =
\bpi_r(l_0+1,\iota_{\boj^+(j)}(l_0+1)-1)_{l_0}.$$ Since $l_0<l$,
the induction hypothesis gives  $$\bpi_r(\boj)_{l_0} =
\bpi_r(\boj^+(j))_{l_0}$$ which implies that
$$\iota_{\boj^+(j)}(l_0+1) = i-r-n+l_0+\iota_\boj(l_0)+2.$$ It
follows that $\bpi_r(\boj^+(j))_l = \bpi_r(\boj)_l$ and the proof
of the Proposition is complete.\hfill\qedsymbol

\subl{}\label{mrep} {\it Proof of Proposition \ref{complete}.}

For this proof it is necessary to indicate that
$\supp^{\pm}(\boj)$  depends on $r$ and  so from now on, if
$\boj\in\bj_{k,r}$ we  denote the set $\supp^\pm(\boj)$ by
$\supp^\pm_r(\boj)$. We also assume  that  the representatives of
the equivalence classes of $\bj_{k,r}$ are chosen so that
$\supp_r^-(\boj) = \emptyset$. We will show that $$\dim
V_{\bpi_r(\boj)\bpi_s(\boj,\pm)} \ge |\bar\boj|= 2^{|\supp^+_r(\boj)|}
\ \ \forall \  \boj\in\bj_{r} \ \ \text{satisfying}\ \
\supp_r^-(\boj)=\emptyset.$$ We proceed by induction on
$M=(i-r)/2$ noting that  induction starts at $M=0$. For the
inductive step, assume that we know the result for all $M'<M$, i.e
we know that for all $s\in I_i$, $s>r$ we have $$\dim
V_{\bpi_s(\boj')\bpi_s(\boj',\pm)} \ge 2^{|\supp^+_s(\boj')|} \
\forall \ \boj'\in\bj_{s} \ \ \text{satisfying}\ \
\supp_s^-(\boj')=\emptyset. $$ Given $j\in \supp^+_r(\boj)$, let
$\boj'\in\bj_{k-1,r+2}$ be the partition whose support is given
by,  $$\supp(\boj') =\{j'+2: r<j'<\sigma^+_\boj(j), \
j'\in\supp(\boj)\}\cup \{j'\in\supp (\boj):
\sigma^+_\boj(j)<j'<n\}.$$  Observe that
$$\supp_{r+2}^-(\boj')=\emptyset,$$ and $$ \supp_{r+2}^+(\boj') =
\{j'+2: r<j'<\sigma^+_\boj(j), \ j'\in\supp^+_r(\boj)\}\cup
\{j'\in\supp^+_r(\boj): \sigma^+_\boj(j)<j'<n\}.$$ This gives,
\begin{align*}
\sigma_\boj^+(j)>j & \implies |\supp_{r+2}^+(\boj')| = |\supp_{r}^+(\boj)|,\\
\sigma_\boj^+(j)=j &\implies  |\supp_{r+2}^+(\boj')| = |\supp_{r}^+(\boj)|-1.
\end{align*}
The inductive step is proved if we show that,
\begin{align*}
\sigma_\boj^+(j)>j & \implies \dim V_{\bpi_r(\boj)\bpi_r(\boj,\pm)}\ge \dim V_{\bpi_{r+2}(\boj')\bpi_{r+2}(\boj',\pm)},\\
\sigma_\boj^+(j)=j &\implies  \dim V_{\bpi_r(\boj)\bpi_r(\boj,\pm)}\ge 2\dim V_{\bpi_{r+2}(\boj')\bpi_{r+2}(\boj',\pm)}.
\end{align*}
Set $l=\sigma^+_\boj(j)$. Using Lemma \ref{wrjDn} we find that,
$$T_{r+2,l}(\bpi_{r+2}(\boj')\bpi_{r+2}(\boj',\pm)) =
\bpi_r(\boj)\bpi_r(\boj,\pm)\balpha_{l,q^{2n-i-l+2\iota_\boj(l)-4}},$$
and that $$(T_{r+2,l}(\bpi_{r+2}(\boj')\bpi_{r+2}(\boj',\pm)))_l =
(1-q^{2n-i-l+2\iota_\boj(l)-4}u)(1-q^{i+l-2r-2\iota_\boj(l)}u).$$
It now follows from Proposition \ref{dim} that $\dim
V_{\bpi_r(\boj)\bpi_r(\boj,\pm)}\ge \dim
V_{\bpi_{r+2}(\boj')\bpi_{r+2}(\boj',\pm)}$. Further,
$$\sigma^+_\boj(j)=j \implies
(T_{r+2,l}(\bpi_{r+2}(\boj')\bpi_{r+2}(\boj',\pm)))_{l} =
(1-q^{n-r-2}u)^2,$$ and we are done by using Proposition \ref{dim}
once more. \hfill\qedsymbol

\bibliographystyle{amsplain}

\begin{thebibliography}{10}

\bibitem{Bo}
N.~Bourbaki, Groupes et alg\`ebres de {L}ie, Chapitres 4,5,6, Hermann, Paris (1968).

\bibitem{BP}
P.~Bouwknegt, K.~Pilch, {\em On deformed W-algebras and quantum affine algebras},
Adv. Theor. Math. Phys. (1998), no.~2, 357-397.

\bibitem{C1} V.~Chari,
{\em Braid group actions and tensor products}, Internat. Math.
Res. Notices (2002), no. 7, 357--382.

\bibitem{CM}
V.~Chari and A.~Moura, {\em Characters and blocks for finite-dimensional representations of quantum affine algebras},
Internat. Math Res. Notices (2005), no. 5, {257--298}.

\bibitem{CPbanff}
V.~Chari and A.~Pressley, {\em Quantum affine algebras and their representations}, in {R}epresentations of groups,
CMS Conf. Proc (Banff, AB, 1994), {\bf 16} (1995), 59--78.

\bibitem{FM}
E. Frenkel\ and\ E. Mukhin, {\em Combinatorics of $q$-characters
of finite-dimensional representations of quantum affine algebras},
Comm. Math. Phys. {\bfseries 216} (2001), no.~1, 23--57

\bibitem{FR}
E.~Frenkel and N.~Reshetikhin, {\em The $q$-characters of
representations of
quantum affine algebras and deformations of $\mathcal{W}$-algebras}, Recent
developments in quantum affine algebras and related topics (Raleigh, NC,
1998), Contemp. Math., 248, Amer. Math. Soc., Providence, RI, 1999,
pp.~163--205.

\bibitem{FH}
W.~Fulton and J.~Harris, Representation Theory - A first course, GTM 129, Springer, (1991).


\bibitem{He}
D.~Hernandez, {\em Monomials of q and q,t-characters for non
simply-laced quantum affinizations}, preprint QA/0404187.

\bibitem{H} J.~E.~Humphreys,
{\em Reflection groups and Coxeter groups}, Cambridge Studies in
Advanced Mathematics {\bf 29}, Cambridge University Press, 1990.



\bibitem{Na} H.~Nakajima, {\em $t$-analogs of $q$-characters of quantum affine algebras of
type $A\sb n,D\sb n$}. Combinatorial and geometric representation
theory (Seoul, 2001), Contemp. Math. {\bf 325} (2003), 141--160.

\end{thebibliography}

\end{document}